\documentclass[a4paper,10pt]{amsart}

\usepackage[all]{xy}
\usepackage{amssymb}
\usepackage{hyperref}
\usepackage{epic}
\usepackage{graphicx}

\newcommand{\A}{{ \rm Aut }}

\newcommand{\F}{{\mathbb{F}}}

\newcommand{\Z}{{\mathbb{Z}}}

\newcommand{\PGL}{{\mathrm{PGL}}}

\newtheorem{pro}{Proposition}[section]
\newtheorem{lemma}[pro]{Lemma}

\newtheorem{cor}[pro]{Corollary}
\newtheorem{theorem}[pro]{Theorem}

\title{On abelian automorphism groups of Mumford curves}
\author{ Aristides Kontogeorgis \and Victor Rotger}

\address{
Department of Mathematics, University of the \AE gean, 83200 Karlovassi, Samos, 
Greece\\ { \texttt{\upshape http://eloris.samos.aegean.gr}}
}
\email{kontogar@aegean.gr}

\address{Departament de Matem\`{a}tica Aplicada IV, Universitat Polit\`{e}cnica de Catalunya, Av. V\'{\i}ctor Balaguer s/n, E-08800 Vilanova i la Geltr\'u, Spain.}
\email{vrotger@ma4.upc.edu}

\subjclass[2000]{11G18, 14G35}

\keywords{Mumford curve, Schottky group, graph, automorphism}

\begin{document}
\bibliographystyle{amsplain}
\date \today
\maketitle

\begin{abstract}
We use rigid analytic uniformization by Schottky groups to give 
a bound for the order of the abelian subgroups of the automorphism 
group of a Mumford curve in terms of its genus. 
\end{abstract}

\section*{Introduction}

Let $X$ be a smooth irreducible projective algebraic curve of genus $g\geq 2$ over a field $k$. 
The automorphism group $\A (X)$ is always finite and it is an interesting problem to determine its size with respect to the genus. When the ground field $k$ has characteristic $0$, it is known that the Hurwitz bound holds:
\begin{equation} \label{hurwitz}
|\A(X)| \leq 84(g-1).
\end{equation}

Moreover, this bound is best possible in the sense that there exist curves of genus $g$ 
that admit $84(g-1)$ automorphisms for infinitely many different values of $g$.

When $\mathrm{char }(k) = p>0$, $|\A(X)|$ is bounded by a polynomial of degree four in $g$. In fact, it holds that
\[
|\A(X)| \leq 16g^4,
\]
\cite{StiI}, provided $X$ is not any of the 
Fermat curves $x^{q+1}+y^{q+1}=1$, $q=p^n$, $n\geq 1$, which have even larger automorphism group \cite{Leopoldt:96}. 

For Mumford curves $X$ over an algebraic extension of the $p$-adic field $\mathbb{Q}_p$,
F. Herrlich \cite{Herrlich:80} was able to improve Hurwitz's bound (\ref{hurwitz}) by showing that  actually 
\[
|\A(X) |\leq 12(g-1),
\] 
provided $p \geq 7$.

Moreover, the first author in a joint work with G. Cornelissen and F. Kato \cite{CKK} proved 
that a bound of the form 
\[
|\A(X)| \leq \max \{12(g-1), 2\sqrt{g} (\sqrt{g}+1)^2 \}
\]
holds for Mumford curves defined over non-archimedean valued fields of characteristic $p>0$.

For ordinary curves $X$ over an algebraically closed field of characteristic $p>0$, Guralnik and Zieve \cite{GZ} announced that there exists a sharp bound of the order of $g^{8/5}$ for $|\A (X)|$. 

In \cite{Nakajima-abel}, S. Nakajima employs the Hasse-Arf theorem to prove that 
\[
|\A(X)| \leq 4g+4
\]
for any algebraic curve $X$ whose group of automorphisms is abelian. 

The results of Herrlich compared to those of Hurwitz and those of \cite{CKK} compared to Guralnik-Zieve's indicate
that if we restrict ourselves to Mumford curves with abelian automorphism group a stronger bound than the one of Nakajima should be expected.
 
The aim  of this note is studying the size of the abelian subgroups of the automorphism group $\A(X)$ of a Mumford curve over a complete field $k$ with respect to a non-archimedean valuation. These curves are rigid analytically uniformized by a Schottky group $\Gamma \subset \PGL_2(k)$ and their automorphism group is determined by the normalizer $N$ of $\Gamma$ in $\PGL_2(k)$.

Our results are based on the {\em Gauss-Bonnet} formula of Karass-Pietrowski-Solitar, 
which relates the rank of the free group $\Gamma$ to the index $[N:\Gamma ]$ and on the 
characterization of the possible abelian stabilizers $N_v\subset N$ of the vertices $v\in \mathcal T_k$ on the Bruhat-Tits tree of $k$ acted  upon by  the group $N$.

\vspace{0.2cm}

{\bf Aknowledgement} The authors would like to thank the referee for his  remarks  and corrections.

\section{Abelian automorphism groups of Mumford Curves}

 For a nice introduction to the theory of automorphism groups of Mumford curves we refer the
reader to \cite{CKNotices}.

Let $k$ be a complete field with respect to a non-archimedean valuation. Let $\bar{k}$ denote the residue field of $k$ and write $p=\mathrm{char }(\bar k)$ for its characteristic. Choose a separable closure $K$ of $k$. 

Let $\Gamma \subset \PGL_2(k)$  be a Schottky group, that is, a discrete finitely generated subgroup consisting entirely of hyperbolic elements acting on $\mathbb{P}^1_k$ with limit set $\mathcal{L}_\Gamma$ (cf.\,\cite{GP}). By a theorem of Ihara, $\Gamma $ is a free group. The rigid analytic curve  
\[
\Gamma \backslash (\mathbb{P}^1_k - \mathcal{L}_\Gamma)
\]
turns out to be the analytic counterpart of a smooth algebraic curve of genus $g = \mathrm{rank} (\Gamma )\geq 1$ over $k$ which we shall denote $X_{\Gamma }/k$. In a fundamental work, D. Mumford \cite{Mumford:72} showed that the stable reduction of $X_{\Gamma }$ is a $\bar k$-split degenerate curve:  all its connected components are rational over $\bar k$ and they meet at ordinary double points rational over $\bar k$. Conversely, he showed that all such curves admit a rigid analytic uniformization by a Schottky subgroup of $\PGL_2(k)$. 

Let $\mathcal{T}_k$ denote the Bruhat-Tits tree of $k$. The set of ends of $\mathcal{T}_k$ is in one-to-one correspondence with the projective line $\mathbb{P}^1(k)$; we thus identify $\mathbb{P}^1(k)$ with the boundary of $\mathcal{T}_k$. 

Let $N$ be a finitely generated discrete subgroup of $\PGL_2(k)$ that contains $\Gamma$ as a normal subgroup of finite index. The group $N$ naturally acts on $\mathcal{T}_k$. By taking an appropriate extension of $k$, we may assume that all fixed points of $N$ in the boundary
are rational. In turn, this implies that $N$ acts on $\mathcal{T}_k$ without inversion.

\begin{theorem}\cite[p.\,216]{GP} \cite{CKK}
The group $G=N/\Gamma$ is a subgroup of the automorphism group of the Mumford curve $X_{\Gamma }$. If $N$ is the normalizer of $\Gamma$ in $\PGL_2(k)$ then $G=\A(X_{\Gamma})$.
\end{theorem}

For every vertex $v$ on $\mathcal{T}_k$ let $N_v$ be the stabilizer of $v$ in $N$, that is,
\[
N_v=\{ g\in N: g(v)=v\}.
\]

Let $\mathrm{star}(v)$ denote the set of edges emanating from the vertex $v$.
It is known that $\mathrm{star}(v)$ is in one to one correspondence with elements in $\mathbb{P}^1(\bar{k})$. 
Since $N_v$ fixes $v$, it acts on $\mathrm{star}(v)$ and describes a natural map  
\begin{equation}\label{rho}
\rho:N_v \rightarrow \PGL_2(\bar{k}).
\end{equation}

See \cite[Lemma 2.7]{CKK}) for details. The kernel of $\rho$ is trivial unless $N_v$ is isomorphic to  the semidirect product of a cyclic group with an elementary abelian 
group. In this case, $\mathrm{ker}\rho$ is an elementary abelian $p$-group \cite[Lemma 2.10]{CKK}.

Assume that $g\geq 2$. This implies that $\Gamma$ has finite index in $N$.
Since $\Gamma$ has finite index in $N$, both groups $N$ and $\Gamma$ share the same set of limit points  $\mathcal{L}$.
We shall denote by $\mathcal{T}_N$ the subtree of  $\mathcal{T}_k$  whose end points are the limit points of $\mathcal{L}$.

The tree $\mathcal{T}_N$ is acted on by $N$ and we can consider the quotient graph $T_N:= N \backslash \mathcal{T}_N$.
The graph $N \backslash \mathcal{T}_N$ is the dual graph of the intersection graph of the special fibre of the  quotient curve 
\[ X_N = G \backslash X_{\Gamma } = N \backslash (\mathbb{P}^1_k - \mathcal{L}).\] 
Notice that $T_N$ is a tree whenever $X_N$ has genus $0$.

The quotient graph $T_N$ can be regarded as a graph of groups as follows: For every vertex $v$ (resp. edge $e$) of $T_N$, consider a
lift $v'$ (resp. $e'$) in $\mathcal{T}_N$ and the corresponding stabilizer $N_{v'}$ 
(resp. $N_{e'}$). We decorate the vertex $v$ (resp. edge $e$) with the stabilizer $N_{v'}$ (resp. $N_{e'}$).

Let $T$ be a maximal tree of $T_N$ and let $T'\subset \mathcal{T}_N$ be a tree of 
representatives of $\mathcal{T}_N$ mod $N$, {\em i.e.}, a lift of $T$ in $\mathcal{T}_N$.
Consider the set $Y$ of lifts of the remaining edges $T_N - T$ in $\mathcal{T}_N$ such that, 
for every $E \in Y$, the origin $o(E)$ lies in $T'$.

The set $Y=\{E_1,\ldots, E_r\}$ is finite. There exist elements $g_i \in N$ such that 
$g_i(t(E_i)) \in T'$, where $t(E_i)$ denotes the terminal vertex of the edge $E_i$ of $Y$. 
Moreover, the elements $g_i$ can be taken from the free group $\Gamma$. 

The elements $g_i$ act by conjugation on the groups $N_{t(E_i)}$ and impose
the relations $g_i N_{t(E_i)} g_i^{-1} =N_{g_i(t(E_i))}$. Denote by  
$M_i:=N_{t(E_i)}$ and $N_i:=N_{g_i(t(E_i))}$.

According to \cite[Lemma 4, p.\,34]{Serre:80}, the group $N$ can be recovered as the group 
generated by 
\[
N:= \langle N_v, g_i\rangle= \langle g_1,\ldots,g_r , K | \,\mathrm{rel}\, K, g_1 M_1 g_1^{-1} =N_1,\ldots,
  g_r M_r g_r^{-1} =N_r \rangle,
\]
where $K$ is the tree product of $T'$. 

Assume that the tree $T'$ of representatives  has $\kappa$ edges and $\kappa+1$ vertices. Let $v_i$ denote the 
order of the stabilizer of the $i$-th vertex and $e_i$ the order of the stabilizer of the $i$-th edge. If $f_i=|M_i|$, we define the {\em volume} of the fundamental domain as
$$
\mu(T_N):=\left(  \sum_{i=1}^r \frac{1}{f_i} + \sum_{i=1}^\kappa \frac{1}{e_i} -
\sum_{i=1}^{\kappa+1} \frac{1}{v_i}
\right).
$$

Notice that when $r=0$, {\em i.e.} the quotient graph $T_N$ is a tree, this definition coincides with the one 
given in \cite{CKK}.

Karrass, Pietrowski and Solitar proved in \cite{KPS} the following {\em discrete Gauss-Bonnet} theorem:

\begin{pro}\label{KPSth} 
Let $N, T_N, g$ be as above. The following equality holds:
\[| N/\Gamma | \cdot \mu(T_{N}) = g-1.\]
\end{pro}

In order to obtain an upper bound for the group of automorphisms with respect to the genus, we aim for a lower bound for $\mu(T_N)$.
Observe that if we restrict the above sum to the maximal tree $T$ of $T_N$, we deduce the following inequality 
\[
\mu(T):=\sum_{i=1}^\kappa \frac{1}{e_i} -
\sum_{i=1}^{\kappa+1} \frac{1}{v_i} \leq \mu(T_N),
\]
where equality is achieved if and only if $T_N$ is a tree, {\em i.e.}, the genus of $X_{N}$ is $0$.
 
In what follows we pursue lower bounds for $\mu(T)$, where $T$ is a maximal tree. These should 
be lower bounds for $\mu(T_N)$ as well.

\begin{lemma} \label{abel-finite}
Let $G$ be a finite abelian subgroup of $\PGL_2(\F_{p^n})$ acting on $\mathbb{P}^1(\mathbb{F}_{p^n})$. 
Let $S$ be the subset of $\mathbb{P}^1(\bar{\F }_p)$ of ramified points of the cover
\[
\mathbb{P}^1 \rightarrow G\backslash \mathbb{P}^1.
\]  

Then, either
\begin{enumerate}
\item 
$G\simeq \Z/n\Z$, where $(n,p)=1$, $S=\{P_1,P_2\}$ and the 
ramification indices are $e(P_1)=e(P_2)=n$, or
\item
$G\simeq D_2 = \Z/2\Z\times \Z/2\Z$, $p\ne 2$ and $S=\{P_1,P_2,P_3\}$ with 
ramification indices $e(P_1)=e(P_2)=e(P_3)=2$, or
\item
$G\simeq E(r)=\Z/p\Z \times \stackrel{(r)}{\cdots }\times  \Z/p\Z$ for some $r\geq 0$ and
$S=\{P\}$, with ramification index $e(P)=p^r$.
\end{enumerate}
\end{lemma}

\begin{proof}
The finite subgroups of $\PGL_2(\F_{p^n})$ were classified by L.\ E.\ Dickson
(cf.\,\cite[II.8.27]{Hu}, \cite{VaMa}, \cite[Theorem 2.9]{CKK}). The list of abelian groups follows by selecting the abelian groups among the 
possible finite subgroups of $\PGL_2(\F_{p^n})$. Notice that the case $E(r) \rtimes \Z/n\Z$, where $(n,p)=1$ and $n \mid p-1$,
is never abelian. Indeed, this is due to the fact that $\Z/n\Z$ acts on $E(r)$ by means of a primitive $n$-th root of unity \cite[cor 1. p.67]{SeLo}. 
The description of the ramification locus $S$ in each case is given in \cite[th 1.]{VaMa}.
\end{proof}

\begin{lemma}\label{ab}
Let $v$ be a vertex of $T_N$. If the finite group $N/\Gamma$ is abelian, then $N_v$ is abelian. Moreover, the map $\rho:N_v \rightarrow \PGL_2(\bar k)$ is injective unless $N_v=E(r_1)$. In this case, 
$\ker(\rho)\simeq E(r_2)$ for some $r_2\leq r_1$.
\end{lemma}

\begin{proof}
The composition 
\[
N_v \subset N \rightarrow N/\Gamma,
\]
is injective, since it is not possible for an element of finite order to be cancelled out by factoring out 
the group $\Gamma$. Hence $N_v$ is a abelian. The possible kernels of $\varrho$ are collected in \cite[Lemma 2.10]{CKK}. 
\end{proof}

Let $v$ be a vertex of $T_N$ decorated by the group $N_v$ and assume that there exist $s\geq 1$ edges in its star, decorated by groups $N_{e_{\nu }}^v \subset N_v$, $\nu =1, ..., s$. We define the {\em curvature} $c(v)$ of $v$ as
\[
c(v):=\frac{1}{2} \sum_{i=1}^s \frac{1}{|N_{e_{\nu }}^v|} -\frac{1}{|N_v|}.
\]

It is obvious that the following formula holds:
\[
\mu(T)=\sum_{v \in \mathrm{Vert}(T)} c(v).
\] 

In what follows, we shall provide lower bounds for the curvature of each vertex. 

We shall call a tree of groups {\em reduced} if $|N_v| > |N_{e_{\nu }}|$ for all vertices $v$ and edges $e_{\nu }\in \mathrm{star }(v)$.
Notice that, if $N_v =N_e$ for a vertex $v$ and an edge $e\in \mathrm{star }(v)$, then the opposite vertex $v'$ of $e$ is decorated by a group 
$N_{v'} \supseteq N_e$. The contribution of $e$ to the tree product is the amalgam $N_v*_{N_e} N_{v'}=N_{v'}$. This means that $e$ can be contracted without altering the tree product. From now on we shall assume that {\em the tree $T$ is reduced}. 

For an element $\gamma\in N$, define the {\em mirror} of $\gamma$ to be
the smallest subtree $M(\gamma)$ of $\mathcal{T}_k$ generated by the point-wise fixed vertices of $\mathcal{T}$ by $\gamma$.

Let $\gamma \in N$ be an elliptic element ({\em i.e.}, an element of $N$ of finite order with 
two distinct eigenvalues of the same valuation). Then $\gamma $ has two fixed points in 
$\mathbb{P}^1(k)$ and $M(\gamma)$ is the geodesic connecting them.

If $\gamma \in N$ is a parabolic element ({\em i.e.}, an element in $N$ having a single 
eigenvalue), then it has a unique fixed point $z$ on the boundary $\mathbb{P}^1(k)$.

\begin{lemma}
Let $P_1,P_2,Q_1,Q_2$ be four distinct points on the boundary of $\mathcal{T}_k$. Let $g(P_1,P_2)$, $g(Q_1,Q_2)$ 
be the  corresponding geodesic on $\mathcal{T}_k$ connecting $P_1,P_2$ and $Q_1,Q_2$ respectively.
For the intersection of the geodesics $g(P_1,P_2)$ and $g(Q_1,Q_2)$ there are the following possibilities:
\begin{enumerate}
\item $g(P_1,P_2)$, $g(Q_1,Q_2)$ have empty intersection.
\item $g(P_1,P_2)$, $g(Q_1,Q_2)$ intersect at only one vertex of $\mathcal{T}_k$.
\item $g(P_1,P_2)$, $g(Q_1,Q_2)$ have a common interval as intersection.
\end{enumerate}
\end{lemma}

\begin{proof}
It immediately follows from the fact that $\mathcal{T}_k$ is simply-connected. 
\end{proof}

We refer to \cite[prop. 3.5.1]{KatoOrb} for a detailed description on the arrangement 
of the geodesics with respect to the valuations of the cross-ratio of the points $P_1,P_2,Q_1,Q_2$.

\begin{lemma} \label{lemma2.6help}
Two non-trivial elliptic elements $\gamma, \gamma'\in\PGL_2(k)$ have the same set of fixed 
points in $\mathbb{P}^1(k)$ if and only if $\langle \gamma, \gamma'\rangle $ is a cyclic group.
\end{lemma}

\begin{proof}
If $\gamma $ and $\gamma'$ generate a cyclic group, there exists an element $\sigma$ such that $\sigma^{i}=\gamma$ and 
$\sigma^{i'}=\gamma'$ for some $i, i'\geq 1$. Since any non-trivial elliptic element has exactly two fixed points, it is immediate that $\gamma$, $\gamma'$ and $\sigma$ have the same set of fixed points.

Conversely if $\gamma,\gamma'$ have the same set of fixed points, say $0,
\infty$, then a simple computation shows that $\gamma,\gamma'$ are of the 
form
$$\gamma=\begin{pmatrix}  a & 0 \\ 0 & d \end{pmatrix} \mbox{ and } \gamma'=\begin{pmatrix} a' & 0 \\ 0 & d' \end{pmatrix},$$
where $a/d$ and $a'/d'$ are roots of unity. Hence there exists $\sigma \in \PGL_2(k)$ such that 
$\sigma^i=\gamma$ and $\sigma^{i'}=\gamma'$.
\end{proof}

\begin{lemma} \label{inter-pos}
Assume that $N/\Gamma $ is an abelian group and let $\gamma,\gamma'\in N$, $\gamma \ne \gamma'$, be elements of prime-to-$p$ finite order. If $M(\gamma)\cap M(\gamma')\ne \emptyset $, then $\langle \gamma, \gamma'\rangle $ is isomorphic to either $D_2$ or a cyclic group.
\end{lemma}

\begin{proof}
By Lemma \ref{ab}, the stabilizers $N_v$ of those vertices $v$ such that $(|N_v|,p)=1$ are abelian 
subgroups of $\PGL_2(\bar k)$. Let $F_{\gamma}$ and $F_{\gamma'}$ denote the sets of the fixed
points of $\gamma$ and $\gamma'$ in $\mathbb{P}^1(k)$, respectively.  

If $M(\gamma )=M(\gamma')$ then $F_{\gamma}= F_{\gamma'}$ and it follows from Lemma \ref{lemma2.6help} that $\langle \gamma, \gamma' \rangle$ is a cyclic group. 

On the other hand, if $M(\gamma )\neq M(\gamma')$, then $F_{\gamma}\neq F_{\gamma'}$ and $\langle \gamma, \gamma' \rangle$  cannot be cyclic, again by Lemma \ref{lemma2.6help}. In this case, any vertex $v\in M(\gamma) \cap M(\gamma')$ is fixed by $\langle \gamma, \gamma' \rangle$, which must be isomorphic to $D_2$ by Lemma \ref{abel-finite}. 
\end{proof}

\begin{lemma} \label{trivial-p-edge}
Assume that $N/\Gamma $ is an abelian group. If $N_v$ is a $p$-group for some vertex $v$ in $T$, then $N_e=\{1\}$ for all $e\in \mathrm{star }(v)$. 
\end{lemma}

\begin{proof}
Assume that $v\in \mathrm{Vert}(T)$ is lifted to $v' \in \mathrm{Vert}(\mathcal{T}_N)$ and that $N_{v'}$ is an 
elementary abelian group. Recall the map $\rho:N_{v'} \rightarrow \PGL_2(\bar{k})$,  which describes the action of $N_{v'}$ on $\mathrm{star}(v')$. 

If $\mathrm{Im}(\rho)=\{\mathrm{Id}_{\PGL_2(\bar{k})}\}$, then every edge $e'\in \mathrm{star}(v')$ would be fixed by $N_{v'}=\mathrm{ker}(\rho)$. 
If one of  the edges $e'\in \mathrm{star}(v')$  were reduced in $T$ to an edge $e$  with  non trivial stabilizer then  it would follow that $N_v=N_e$ and the tree would not be reduced.

Suppose now that  $\mathrm{Im}(\rho)\neq\{\mathrm{Id}_{\PGL_2(\bar{k})}\}$. 
By lemma \ref{abel-finite} there exists exactly  one edge $e'$ in the star of $v'$ which 
is fixed by the whole group $\mathrm{Im}(\rho)$ and all other edges emanating from $v'$ are not fixed by $\mathrm{Im}(\rho)$. 
Therefore the edge $e'$ is fixed by the whole group $N_{v'}$. If the edge $e'$ were  reduced in $T$ to an edge $e$ with non trivial stabilizer, then 
it would follow that $N_v=N_e$ and the tree would not be reduced.

Assume now that there exist two vertices $v_1,v_2$ on $T$ joint by an edge $e$ such that $N_{v_1},N_{v_2}$ are 
elementary abelian groups and $N_e$ is a nontrivial proper subgroup both of $N_{v_1}$ and $N_{v_2}$. Let us show that this 
can not happen.

Let $\sigma, \tau$ be two commuting parabolic elements of $\PGL_2(\bar{k})$. They 
fix a common point in the boundary of $\mathbb{P}^1(k)$. Indeed, every parabolic element fixes a single point in the boundary.
If $P$ is the unique fixed point of $\sigma$, then 
\[\sigma(\tau P)=\tau (\sigma P)=\tau P\]
and $\tau(P)$ is fixed also by $\sigma$. Since the fixed point of $\sigma$ in the boundary is unique, we have $\tau(P)=P$.

Let $v_1',v_2'$ be two lifts of $v_1,v_2$ on the Bruhat-Tits tree. The  apartment $[v_1',v_2']$  is contracted to the edge $e$ and it 
is fixed by $N_e$, but not by a larger subgroup.
 
Since $N_e$ is contained in both abelian groups $N_{v_1'},N_{v_2'}$, all parabolic elements in $N_{v_1'},N_{v_2'}$
have the same fixed point $P$ in the boundary $\mathbb{P}^1(k)$. Therefore, the apartment $[v_1',P[$ 
(resp. $[v_2',P[$) is fixed by $N_{v_1'}$ (resp. $N_{v_2'}$). Moreover, the 
apartments $[v_1',P[$,$[v_2',P[$ have nonempty intersection. Since the Bruhat-Tits tree is simply connected, the apartment  $[v_1',v_2']$ intersects $[v_2',P[ \cap [v_1',P[$ at a bifurcation point $Q$:
\[[v_1',v_2'] \cap \big([v_2',P[ \cap [v_1',P[ \big)=\{Q \}.\]

The  point $Q$ is then  
fixed by $N_{v_1'}$ and $N_{v_2'}$ and it is on the apartment $[v_1',v_2']$, a contradiction.
\end{proof}

\begin{lemma} \label{cont-estim}
Let $v$ be a vertex in $T$. If $c(v)>0$, then $c(v)\geq \frac{1}{6}$. Let $s = \# \mathrm{star }(v)$ and let $N_{e_\nu}^v$ denote the stabilizers of the edges in the star of $v$ for $\nu=1,\ldots,s$. It holds that $c(v)=0$ if and only if: 
\begin{enumerate}
\item $N_v=D_2$, $s=1$, $N_{e_1}^v=\mathbb{Z}_2$ or
\item $N_v=\mathbb{Z}_2$, $s=1$, $|N_{e_1}^v|=1$.
\end{enumerate}
We have
\begin{itemize}
\item $c(v)=\frac{1}{6}$ if and only if $N_v=\mathbb{Z}_3$ and $s=1$,
\item $c(v)=\frac{1}{4}$ if and only if 
$N_v=D_2$ with $s=2$ and $|N_{e_1}^v|=|N_{e_2}^v|=|2|$, or $N_v= D_2$ with $s=1$ and 
$|N_{e_1}^v|=1$.
\end{itemize}
 
In the remaining cases we have $c(v)\geq \frac{1}{3} $.
\end{lemma}

\begin{proof}
$\bullet$  Assume that $N_v=D_2$. Then  $c(v)\geq 0$. Equality $c(v)=0$ holds only when $s=1$ and the only edge leaving $v$ is decorated by a group of order $2$.
If we assume that $c(v) > 0$, then 
\[
\frac{1}{4} \leq c(v)
\]
and equality is achieved if $s=2$ and $|N_{e_1}|=|N_{e_2}|=2$, or if $s=1$ and $|N_{e_1}^v|=1$.

$\bullet$ Assume that $N_v=\mathbb{Z}_n$. Then  
\[
 c(v)= \sum_{i=\nu}^s  \frac{1}{2|N_{e_\nu}|} - \frac{1}{n}.
\]

By Lemma \ref{inter-pos}, the stabilizer of each edge in the star of $v$ is trivial.
Indeed, if there were an edge $e \in \mathrm{star(v)}$ with $N_e >\{1\}$, then $e$ would be fixed by a cyclic group $\mathbb{Z}_m$, where $m \mid  n$.  
Let $\sigma$ be a generator of $\mathbb{Z}_n$ and let $\sigma^\kappa$ be the generator of $\Z_m$.
The elements $\sigma,\sigma^\kappa$ have the same fixed points. Hence a lift of the edge $e$ in $\mathcal{T}_N$ would 
lie on the mirror of $\sigma$. But then $N_e=N_v$ and this is not possible by the  
reducibility assumption. See also \cite[Lemma 1]{Herrlich:80}.

If $n > 2$, then 
\[
\frac{1}{6} \leq \frac{n-2}{2n}\leq \frac{sn-2}{2n}=s\frac{1}{2}-\frac{1}{n} \leq c(v),
\]
and equality holds only if $s=1,n=3$.
If $n=2$ and $c(v)>0$ then $s\geq 2$, and $c(v)=s\frac{1}{2}-\frac{1}{2} \geq \frac{1}{2}$.

$\bullet$ 
Assume that $N_v=E(r)$. Then $s=1$ and it follows from Lemma \ref{trivial-p-edge} that $|N_{e_1}|=1$. If $p^r=2$ then 
$c(v)=0$. Hence if $c(v)>0$ then  $p^r>2$ and 
\[
\frac{1}{6} \leq c(v)=\frac{1}{2}-\frac{1}{p^r}=\frac{p^r-2}{2p^r},
\]
and equality holds only if $p^r=3$.
\end{proof}

\begin{theorem} \label{mainth} Assume that $N/\Gamma$ is abelian. If $N$ is  neither  isomorphic to $\mathbb{Z}_2 * \mathbb{Z}_3$ nor $D_2* \mathbb{Z}_3$, then 
\[ 
|N/\Gamma| \leq 4 (g-1).
\]
If we exclude the groups of Table \ref{table} then 
\[
|N/\Gamma| \leq 3 (g-1).
\] 
The case $N=\mathbb{Z}_2 * \mathbb{Z}_3$ gives rise to a curve of genus $2$ whose automorphism group is 
a cyclic group of order $6$. The case 
$N=D_2* \mathbb{Z}_3$ gives rise to a curve of genus $3$ with automorphism group $D_2 \times \mathbb{Z}_3$.\begin{table}
\[
\begin{array}{c|c|c}
N & N/\Gamma & g \\
\hline
D_2 * \mathbb{Z}_4  &   D_2 \times  \mathbb{Z}_4   &  2  \\
\mathbb{Z}_2 * \mathbb{Z}_4  &   \mathbb{Z}_2 \times  \mathbb{Z}_4   & 2   \\
D_2 * D_2  &   \mathbb{Z}_2 ^4   &   4 \\
\mathbb{Z}_2      * D_2  &   \mathbb{Z}_2 ^3   &  2  \\
        D_2 *_{\mathbb{Z}_2}    D_2 * _{\mathbb{Z}_2 } D_2 &  \mathbb{Z}_2 ^4  & 2
\end{array}
\]
\caption{ \label{table} }
\end{table}
\end{theorem}

\begin{proof}
Since $g\geq 2$ and therefore $\mu(T_N) > 0$, we have by Proposition \ref{KPSth} that 
\[
|\A(X)| =\frac{1}{\mu(T_N)} (g-1) \leq  \frac{1}{\mu(T)} (g-1)  \leq \frac{6}{\lambda} (g-1),
\]
where 
\[
\lambda=\#\{ v \in \mathrm{Vert}(T_N): c(v)>0 \}.
\]

If $\lambda\geq 2$ the result follows. Assume that there is only one vertex $v$ such that $c(v)>0$.
Since $g\geq 2$, there exist other vertices $v'$ on the tree $T_N$ but their contribution is $c(v')=0$.
Notice that if we contract  a tree along an edge connecting the vertices $v_1,v_2$
forming a new vertex $v$ then $c(v_1)+c(v_2)=c(v)$. Therefore, one can check 
that $c(v) \geq 0$, using lemma \ref{abel-finite}.

{\bf Case 1:}  $c(v)=\frac{1}{6}$. Then $N_v=\mathbb{Z}_3$ and there exists a single edge $e$ at the star of $v$. Let $v'$ denote the terminal vertex of $e$. Since $c(v')=0$ if and only if there exists a single edge leaving $v'$, the only 
possibilities for $N$ are $N=D_2 * \mathbb{Z}_3$ and $N=\mathbb{Z}_2 * \mathbb{Z}_3$.
Since $N/\Gamma $ is abelian, the group $\Gamma_1:=[D_2,\mathbb{Z}_3]$ (resp. 
$\Gamma_1:=[\mathbb{Z}_2,\mathbb{Z}_3])$ is contained in $\Gamma$. 
According to \cite[Lemma 6.6]{CKK}, $\Gamma_1$ is a  maximal free subgroup of $N$ and thus 
$\Gamma=\Gamma_1$. The rank of $\Gamma$ is $(4-1)(2-1)=3$ in the first case and $(3-1)(2-1)=2$ in the second. Therefore 
the first amalgam  gives rise to a curve of genus $3$ with automorphism group $D_2 \times \mathbb{Z}_3$
and the second gives rise to a curve of genus $2$ with automorphism group $\mathbb{Z}_2 \times 
\mathbb{Z}_3 \cong \mathbb{Z}_6$.

{\bf Case 2:} $c(v)=\frac{1}{4}$. This occurs only when $N_v=\mathbb{Z}_4,D_2$ and $s=1$ 
or $N_v=D_2$, $s=2$, $|N_{e_1}|=|N_{e_2}|=\frac{1}{2}$. The possible groups are given in Table:\ref{table}.

In this case we have the following bound
\[
|\A(X)| \leq \frac{1}{\mu(T)} (g-1)  \leq 4(g-1).
\]

{\bf Case 3:} $c(v)\geq \frac{1}{3}$. Similarly as above we obtain that 
\[
|\A(X)| \leq \frac{1}{\mu(T)} (g-1)  \leq \frac{3}{\lambda} (g-1) \leq 3 (g-1).
\]

\end{proof}

{\bf Example:} {\em Subrao Curves.}

Let $(k,|\cdot|)$ be a complete field of characteristic $p$ with respect 
to a non-archimedean norm $|\cdot|$. Assume $\mathbb{F}_q \subset k$, for some
$q=p^r$, $r\geq 1$. Define the curve:
\[
(y^q-y)(x^q-x)=c,
\]
with $|c|<1$.  This curve was introduced by Subrao  in \cite{subrao} and 
it has a large automorphism group compared to the genus.
This curve is a Mumford Curve \cite[p.\,9]{CKK} and has {\em chessboard} reduction \cite[par. 9]{CKK}.
It is a curve of genus $(q-1)^2$ and admits the group $G:=\mathbb{Z}_p^r \times \mathbb{Z}_p^r$ as a subgroup of the automorphism group. 
The group $G$ consists of  the automorphisms $\sigma_{a,b}(x,y)=(x+a,y+b)$ where $(a,b)\in \mathbb{F}_q\times \mathbb{F}_q$. The discrete group $N'$ 
corresponding to $G$ is given by $\mathbb{Z}_{p^r} * \mathbb{Z}_{p^r}$ and the free subgroup $\Gamma$ giving the Mumford uniformization is 
given by the commutator $[ \mathbb{Z}_{p^r} , \mathbb{Z}_{p^r}]$, which is of rank $(q-1)^2$ \cite{Serre:80}. 

Our bound is given by 
\[
q^2 =|G| \leq 2 (g-1)= 2 (q^2-2q).
\] 
Notice that the group $N'$ is a proper subgroup of  the normalizer of $\Gamma$ in $\PGL_2(k)$, since the full automorphism group of the curve is 
isomorphic to $\mathbb{Z}_p^{2r} \rtimes D_{p^r-1}$ \cite{CKK}. $\Box $

\subsection{Elementary abelian groups}

\begin{pro}\label{mainpro}
Let $\ell $ be a prime number and let $X_{\Gamma }/k$ be a Mumford curve over a non-archimedean local field $k$ such that $p=\mathrm{char}(\bar k)\ne \ell $.  Let $A \subset \A(X_{\Gamma })$ be a subgroup of the group of automorphisms of $X_{\Gamma }$ such that $A\simeq \Z/\ell \Z \times \Z/\ell \Z\times  \cdots \times \Z/\ell \Z$. 

If $\ell =2$ then all stabilizers of vertices and edges of the quotient graph $T_N$ are subgroups of $\Z/2\Z\times \Z /2\Z$ and $\mu(T_N)=a/4$ for some $a\in \mathbb{Z}$. If $\ell >2$ then all stabilizers of vertices and edges of $T_N$ are subgroups of $\Z/\ell \Z$ and $\mu(T_N)=a/\ell $ for some $a\in \mathbb{Z}$.
\end{pro}

\begin{proof}
Let $A \subset \A(X_{\Gamma})$. There is a discrete finitely generated subgroup $N'\subset N$ such that $\Gamma 
\lhd N'$ and $N'/\Gamma=A$. Let $N_v$ be the stabilizer of a vertex in $T_N$ and let $N'_v=N_v \cap N'$. 
The composition 
\[
N_v' \subset N_v \subset N \rightarrow N/\Gamma,
\]
is injective, since it is not possible for an element of finite order to be 
cancelled out by factoring out 
the group $\Gamma$. 

The map $\rho:N_v' \rightarrow \PGL_2(\bar{k})=\PGL_2(\F_{p^m})$
of Lemma \ref{ab} is injective since $(|N_v'|,p)=1$ and hence we can regard 
$N_v'$ as a finite subgroup of $\PGL_2(\bar{k})=\PGL_2(\F_{p^m})$. 

Assume first that $\ell =2$. Then by Lemma \ref{abel-finite} the only abelian finite subgroups of $\PGL_2(\bar{k})$ for $p\neq 2$ are
$\Z/2\Z$ and the dihedral group of order $4$. Hence $N_v'$ is a 
subgroup of $\Z/2\Z\times \Z /2\Z$.

Since the group $N$ acts without inversions, the stabilizer of a vertex is the 
intersection of the stabilizers of the limiting edges. It again follows that $N_e\subseteq \Z/2\Z\times \Z /2\Z$. 
Finally, we obtain from its very definition that $\mu(T_N)=a/4$ for some $a\in 
\mathbb{Z}$. 

For the case $\ell >2$ we observe that $\Z/\ell \Z$ is the only abelian subgroup of 
$\PGL_2(\F_{p^m})$ and 
it follows similarly that $\mu(T_N)=a/\ell $ for some $a\in \mathbb{Z}$. 
\end{proof}

\vspace{0.2cm}
As an immediate corollary of Proposition \ref{mainpro} we obtain the following formula for the $\ell $-elementary subgroups of the group of automorphisms of Mumford curves.

Notice that the result below actually holds for arbitrary algebraic curves -as it can be proved by applying the Riemann-Hurwitz formula to the covering $X\longrightarrow X/A$.

\begin{cor} \label{main}
Let $\ell \ne \mathrm{char}(\bar{k})$ be a prime number and let $X/k$ be a Mumford curve of genus $g\ge 2$ over a non-archimedean local field $k$. Let $A\subseteq \A(X)$ be a subgroup of the group of automorphisms of $X$ such that $A\simeq \oplus_{i=1}^s \Z/\ell \Z$ for some $s\ge 2$. 
\begin{enumerate}
\item[(i)] If $\ell \ne 2$ then $\ell ^{s-1} \mid g-1$. 
\item[(ii)] If $\ell =2$ then $2^{s-2} \mid g-1$. 
\end{enumerate} 
\end{cor}


\begin{thebibliography}{9}

\bibitem{CKNotices}
G.\ Cornelissen, F.\ Kato, {The {$p$}-adic icosahedron},
{\em Notices Amer.\ Math.\ Soc.\ } \textbf{52:7} (2005), 720-727. 


\bibitem{CKK} 
G.~Cornelissen, F.~Kato, A.~Kontogeorgis, {Discontinuous groups in positive characteristic 
and automorphisms of {M}umford curves}, {\em Math.\ Ann.\ } {\bf 320} (2001), 55-85.

                
\bibitem{GP}    
L.\ Gerritzen, M.\ van der Put, {\em Schottky groups and Mumford curves}, Lect.\ Notes Math.\ {\bf 817}, Springer 1980.
                

\bibitem{GZ}
R.\ M.\ Guralnick,  M.\ E.\ Zieve,  Work on  {\em Automorphisms of Ordinary Curves}, in preparation 
(Talk in Leiden, Workshop on Automorphism of Curves, 18-08-2004).


\bibitem{Herrlich:80}
F.\ Herrlich, {Die {O}rdnung der {A}utomorphismengruppe einer
{$p$}-adischen {S}chottkykurve}, {\em Math.\ Ann.\ } \textbf{246:2} (1979/80), 
125-130. 

\bibitem{Hu}
B.\ Huppert, Endliche Gruppen I, in {\em Die Grundlehren der Mathematischen Wissenschaften,} {\bf 134} 
Springer (1967), 20-25.


\bibitem{KatoOrb} F.\ Kato, Non-archimedean orbifolds covered by Mumford curves, {\em J.\ Algebraic Geom.\ }
{\bf 14:1} (2005), 1-34

\bibitem{KPS}
A.~Karrass, A.~Pietrowski, D.~Solitar, {Finite and infinite cyclic
extensions of free groups}, {\em J.\ Austral.\ Math.\ Soc.\ } \textbf{16} (1973), 458-466. 



\bibitem{Leopoldt:96}
H.\ W.\ Leopoldt, {\"{U}ber die {A}utomorphismengruppe des {F}ermatk\"orpers}, {\em J.\ Number Theory} \textbf{56:2} (1996), 256-282.

\bibitem{Mumford:72}
D.~Mumford, {An analytic construction of degenerating curves over complete local rings}, {\em Compositio Math.\ } {\bf 24} (1972), 129-174.

\bibitem{Nakajima-abel}
S.~Nakajima, {On abelian automorphism groups of algebraic curves},  {\em J.\ London Math.\ Soc.\ }, 
{\bf 36}, (1987), 23-32.

\bibitem{SeLo}
J.-P.\ Serre, \emph{Local Fields}, Grad.\ Texts Math.\ {\bf 67}, Springer 1979.

\bibitem{Serre:80}
J.-P. Serre, \emph{Trees}, Springer 1980.

\bibitem{StiI}
H.\ Stichtenoth, {\"{U}ber die {A}utomorphismengruppe eines
algebraischen {F}unktionenk\"orpers von {P}rimzahlcharakteristik. {I}. {E}ine
{A}bsch\"atzung der {O}rdnung der {A}utomorphismengruppe}, {\em Arch.\ Math.\ }
(Basel) \textbf{24} (1973), 527-544. 

\bibitem{subrao}
D.\ Subrao, {The {$p$}-rank of {A}rtin-{S}chreier curves},
{\em Manuscripta Math.\ } \textbf{16:2} (1975), 169-193. 

\bibitem{VaMa}
R.\ C.\ Valentini, M.\ L.\ Madan,  {A {H}auptsatz of L.\ E.\ Dickson and Artin-Schreier extensions}, {\em  J.\ Reine Angew.\ Math.\ } {\bf 318} (1980), 156-177.

\end{thebibliography}
\end{document}